\documentclass[12pt]{amsart}

\usepackage{amssymb}

\usepackage{enumitem}

\usepackage{graphicx}

\makeatletter
\@namedef{subjclassname@2020}{%
  \textup{2020} Mathematics Subject Classification}
\makeatother

\usepackage[T1]{fontenc}

\newtheorem{theorem}{Theorem}[section]
\newtheorem{corollary}[theorem]{Corollary}
\newtheorem{lemma}[theorem]{Lemma}

\theoremstyle{definition}
\newtheorem{definition}[theorem]{Definition}
\newtheorem{remark}[theorem]{Remark}
\newtheorem{example}[theorem]{Example}
\newtheorem{question}[theorem]{Question}
\newtheorem{cl}{Claim}[section]

\newcommand{\bques}{\begin{question}}
\newcommand{\eques}{\end{question}}
\newcommand{\beg}{\begin{example}}
\newcommand{\eeg}{\end{example}}
\newcommand{\bpf}{\begin{proof}}
\newcommand{\epf}{\end{proof}}
\newcommand{\bcl}{\begin{cl}}
\newcommand{\ecl}{\end{cl}}

\numberwithin{equation}{section}

\newcommand{\diam}{{\operatorname{diam}}}
\newcommand{\id}{{\operatorname{id}}}
\newcommand{\dist}{{\operatorname{dist}}}

\def\be{\begin{equation}}
\def\ee{\end{equation}}

\begin{document}


\baselineskip=17pt


\title[Gromov hyperbolicity in the free quasiworld. I]{Gromov hyperbolicity in the free quasiworld. I}

\author[Q. Zhou]{Qingshan Zhou}
\address{School of Mathematics and Big Data,\\ Foshan University, \\
Foshan, Guangdong 528000,\\
People's Republic of China} 
\email{qszhou1989@163.com; q476308142@qq.com}

\author[S. Ponnusamy]{Saminathan Ponnusamy}
\address{Department of Mathematics,\\ Indian Institute of Technology Madras,\\
Chennai 600036, India}
\address{Department of Mathematics,\\ Petrozavodsk State University, \\
ul., Lenina 33,\\
185910 Petrozavodsk, Russia}
\email{samy@iitm.ac.in}

\date{}

\begin{abstract}
With the aid of a Gromov hyperbolic characterization of uniform domains, we first give an affirmative answer to an open question arisen by V\"ais\"al\"a under weaker assumption. Next, we show that the three-point condition introduced by V\"ais\"al\"a is necessary to obtain quasisymmetry for quasim\"obius maps between bounded connected spaces in a quantitative way. Based on these two results, we investigate the boundary behavior of freely quasiconformal and quasihyperbolic mappings on uniform domains of Banach spaces and partially answer another question raised by V\"ais\"al\"a in different ways.
\end{abstract}

\subjclass[2020]{Primary 30C65, 30F45, 30L10; Secondary 30C20}

\keywords{Quasisymmetric map, quasim\"{o}bius map, bilipschitz map, Gromov hyperbolic domain, uniform domain, freely quasiconformal map, quasihyperbolic metric}

\maketitle

\section{Introduction and main results}\label{sec-1}

Around 1990, V\"ais\"al\"a began a systematic study of quasiconformal and related maps in arbitrary Banach spaces, see \cite{Vai-1,Vai-2,Vai-3, Vai-4,Vai-5}. By invoking the quasihyperbolic metric, a useful tool in geometric function theory introduced by Gehring and his students \cite{GO, GP}, V\"ais\"al\"a established his {\it freely quasiconformal mapping} (FQC) theory by elementary arguments only, since a number of methods and tools in plane or higher dimensional quasiconformal theory are not available in an infinite-dimensional Banach space such as conformal capacity or the moduli of path families, and Whitney decomposition. More importantly, there are several open questions proposed by V\"ais\"al\"a in the free quasiworld to investigate whether or not several well-known results in $\mathbb{R}^n$ can be formulated in a dimension-free manner. These questions have been largely motivated and considered by Huang, Li, Rasila, Vourinen and Wang, see \cite{HLVW,LVX,rt,rt2,ZR} etc. Recently, the relationship between the quasihyperbolic geometry and weakly quasisymmetric maps has been investigated by Huang and Liu \cite{HL} and Huang et al. \cite{HLL}, in quasiconvex metric spaces. More recently, Huang et al. \cite{HRWZ} proved the equivalence of freely quasiconformal mappings and locally quasisymmetric maps in the setting of metric spaces.

After the initiation of freely quasiconformal mapping theory by V\"ais\"al\"a, Heinonen (in the Mathematical Review--MR1709971 of the article of V\"ais\"al\"a \cite{Vai-5}) pointed out that it would be interesting to investigate the relationship between the infinite-dimensional freely quasiconformal mapping theory and the work of Gromov. In \cite{Gr87}, Gromov observed that the asymptotic geometric properties of $n$-hyperbolic space $\mathbb{H}^n$ can be characterized by a simple condition for quadruples of points.  This condition is now used to define a class of metric spaces that is named  Gromov hyperbolic spaces. It forms a large and much studied class of metric spaces and plays an important role in geometric group theory and geometric function theory; see \cite{BB03,BHK,Her06,HSX,K,Vai} and the references therein. More recently, Butler pointed out it is an interesting question whether the theory of uniformizing Gromov hyperbolic spaces could be applied to the free quasiworld, see \cite[Remark 1.8]{Bu}.

In view of these considerations, we study the Gromov hyperbolicity of the quasihyperbolic metric with its applications in V\"ais\"al\"a's free quasiworld. As the first motivation of this paper, we give an affirmative answer to an open question proposed by V\"ais\"al\"a under weaker assumption by virtue of a Gromov hyperbolic characterization of uniform domains. Let us briefly recall the history of this question. In this section, we assume that $E$ and $E'$ are real Banach spaces with dimension at least $2$, and that $D\subsetneq E$ and $D'\subsetneq E'$ are proper domains. We begin with the following equivalent conditions of freely quasiconformal mappings on uniform domains of Banach spaces. In \cite{Vai-2}, V\"ais\"al\"a proved

\begin{theorem}\label{Thm-A}
$($\cite[Theorem 7.18]{Vai-2}$)$
Suppose that $D$ is a $c$-uniform domain and that $f:\,D\to D'$ is freely $\varphi$-quasiconformal. Then the following conditions are quantitatively equivalent:
\begin{enumerate}
  \item $D'$ is $c_1$-uniform,
  \item $f$ is $\theta$-quasim\"obius.
\end{enumerate}
\end{theorem}

The terminology used in Theorem \ref{Thm-A} and in  the rest of this section will be explained in Section \ref{sec-2}. The class of freely quasiconformal mappings  was introduced by  V\"ais\"al\"a in \cite{Vai-1}, see Definition \ref{def-fqc}. In $\mathbb{R}^n$, a homeomorphism between two proper domains is freely $\varphi$-quasiconformal if and only if it is $K$-quasiconformal, where $\varphi$ and $K$ depend only on each other and $n$. In order to establish the theory of quasiconformal mappings in infinite dimensional Banach spaces, V\"ais\"al\"a consider this class of mappings instead of quasiconformal mappings.

The proof of Theorem \ref{Thm-A} depends almost entirely on the following result:

\begin{theorem}\label{Thm-B} $($\cite[Theorem 7.9]{Vai-2}$)$
Suppose that $D\neq E$ and $D'\neq E'$ are $c$-uniform domains and that $f:\,D\to D'$ is $C$-coarsely $M$-quasihyperbolic. Then $f$ extends to a homeomorphism $\overline{f}:\,\overline{D}\to \overline{D'}$ and $\overline{f}$ is $\theta$-quasim\"obius relative to $\partial D$ with $\theta$ depending only on $M,C$, and $c$. In particular, $f|_{\partial D}$ is $\theta$-quasim\"obius.
\end{theorem}

We briefly review the relations between quasiconformal mappings and coarsely quasihyperbolic mappings. It was proved by Gehring and Osgood in \cite[Theorem 3]{GO} that a quasiconformal mapping $f$ between two proper domains in $\mathbb{R}^n$ and its inverse are uniformly continuous with respect to the quasihyperbolic metrics. Because the quasihyperbolic metric is a length metric, we easily see that $f$ is a coarsely quasihyperbolic mapping, i.e., a coarsely bilipschitz mapping  with respect to the quasihyperbolic metrics. Similarly, freely quasiconformal mappings are also coarsely quasihyperbolic.

Subsequently, for the converse of Theorem \ref{Thm-B}, V\"ais\"al\"a raised the following:

\begin{question}\label{q-1} $($\cite[Question 7.19]{Vai-2}$)$ Suppose that $D$ is a $c$-uniform domain and that $f:\,D\to D'$ is $C$-coarsely $M$-quasihyperbolic. If $f$ extends to a homeomorphism $\overline{f}:\,\overline{D}\to \overline{D'}$ and $\overline{f}$ is $\theta$-quasim\"obius relative to $\partial D$, is $D'$   $c'$-uniform ?
\end{question}

Recently, Huang et al. \cite{HLVW} answered Question \ref{q-1} in the affirmative by using the method associated to the quasihyperbolic geometry and the notion of {\it coarse length} of an arc introduced by V\"ais\"al\"a \cite[Section 4]{Vai-2}.

\begin{theorem}\label{Thm-C} $($\cite[Theorem 1.1]{HLVW}$)$
Suppose that $D$ is a $c$-uniform domain and that $f:\,D\to D'$ is $C$-coarsely $M$-quasihyperbolic. Then $f$ extends to a homeomorphism $\overline{f}:\,\overline{D}\to \overline{D'}$ and $\overline{f}$ is $\theta$-quasim\"obius relative to $\partial D$ if and only if $D'$ is $c'$-uniform, quantitatively.
\end{theorem}

Moreover, it is very natural to consider the following.

\begin{question}\label{q-2} Suppose that $D$ is a $c$-uniform domain and that $f:\,D\to D'$ is $C$-coarsely $M$-quasihyperbolic. If $f$ extends to a homeomorphism $\overline{f}:\,\overline{D}\to \overline{D'}$ and $\overline{f}$ is $\theta$-quasim\"obius on $\partial D$, is $D'$   $c'$-uniform ?
\end{question}

Note that each quasim\"obius mapping relative to the boundary is quasim\"obius on the boundary. Therefore, we know that the assumptions mentioned in Question  \ref{q-2} is weaker than that of Question  \ref{q-1}.

Our first main result (Theorem \ref{main-thm-1} below) answers affirmatively Question \ref{q-2} which is also an improvement of Theorem \ref{Thm-C}. We also remark that our proof is based on the Gromov hyperbolic characterization of uniform domains established by V\"ais\"al\"a \cite{Vai}, which states that a domain in a Banach space is uniform if and only if it is Gromov hyperbolic in the quasihyperbolic metric and the norm boundary is naturally quasim\"obius equivalent to the Gromov boundary equipped with a visual metric.

\begin{theorem}\label{main-thm-1}Let $M,c\geq 1$ and $C\geq 0$. Suppose that $D$ is a $c$-uniform domain and that $f:\,D\to D'$ is $C$-coarsely $M$-quasihyperbolic. Then $f$ extends to a homeomorphism $\overline{f}:\,\overline{D}\to \overline{D'}$ and $\overline{f}$ is $\theta$-quasim\"obius on $\partial D$ if and only if $D'$ is $c'$-uniform, where $\theta$ and $c'$ depend only on each other and, $M$, $C$, and $c$.
\end{theorem}

Secondly, we investigate a sufficient condition for quasim\"obius maps to be quasisymmetric, which was introduced by V\"ais\"al\"a in \cite{Vai-0}. Recall that:

\begin{definition}
A homeomorphism $f:\, (X,|\cdot|)\to (Y,|\cdot|)$ between bounded metric spaces is said to satisfy the {\it $\lambda$-three-point condition} for some quadruple $(z_1,z_2,z_3;\lambda)$, if there is a  constant $\lambda\geq 1$ and some points $z_1,$ $z_2,$ $z_3$ in $X$ such that
\[ 
|z_i-z_j|\geq \frac{1}{\lambda}\diam (X)\;\;\mbox{and}\;\; |f(z_i)-f(z_j)|\geq \frac{1}{\lambda}\diam (Y)
\]
for all $i, j\in\{1,2,3\}$ with $i\neq j$.
\end{definition}

It is clear that a homeomorphism  between bounded metric spaces has the $\lambda$-three-point condition for some $\lambda\geq 1$ depending on the locations of these points and the diameters of the spaces. Under this condition, it was shown by V\"{a}is\"{a}l\"{a} in \cite{Vai-0} that every quasim\"{o}bius map between two bounded metric spaces is quasisymmetric in a quantitative way. That is, the control function for quasisymmetric maps depends only on the quasim\"obius control function and the constant $\lambda$ associated with the three-point condition defined as above.

We also give a short proof of this fact by using a useful notion of deformed cross-ratio introduced by Bonk and Kleiner in \cite{BK1}. Recently, aided by this tool, Wang and Zhou \cite{WZ17} systematically studied the relationship between weakly quasim\"obius mappings and quasim\"obius mappings in doubling metric spaces.
Our second issue is to prove that the three-point condition is actually necessary for a quasim\"obius map to be quasisymmetric in a quantitative way.

\begin{theorem}\label{main-thm-2}
Suppose that $f:\,X\to Y$ is a homeomorphism between bounded connected metric spaces. Then $f$ is $\eta$-quasisymmtric if and only if $f$ is $\theta$-quasim\"{o}bius and satisfies the $\lambda$-three-point condition,  where $\eta$ and $\theta$, $\lambda$ depend only on each other.
\end{theorem}

It turns out that Theorem \ref{main-thm-2}, together with Theorem \ref{main-thm-1}, is very helpful to understand the boundary behavior of freely quasiconformal mappings on uniform domains of Banach spaces. Many related and interesting conclusions are given in subsection \ref{sub-2}. For instance, we prove that if two proper domains of Banach spaces are quasim\"obius equivalent with quasisymmetric boundary correspondence, then they are quasisymmetric to each other, see Lemma \ref{s-7}. Moreover, it is shown in Theorem \ref{s-9} that if a domain is freely quasiconformal to a uniform domain with their boundaries quasisymmetrically equivalent, then they are quasisymmetric to each other. In particular, the target domain is uniform. We also remark that this result is an improvement of \cite[Theorem 1.2]{LVX} under weaker hypotheses.

As applications, we also investigate boundary behavior of quasihyperbolic maps on uniform domains and consider an open question raised by V\"ais\"al\"a.

\begin{question}\label{q-3}$($\cite[Question 13]{Vai-5}$)$ Suppose that $D$ and $D'$ are bounded domains with connected boundaries, suppose also that $f:\,D\to D'$ is $L$-quasihyperbolic and that $f$ extends to a homeomorphism $f:\,\overline{D}\to \overline{D'}$ such that $f|_{\partial D}$ is $L$-bilipschitz. Is $f$ $L'$-bilipschitz with $L'=L'(L)$?
\end{question}

Recall that a quasihyperbolic map between proper domains means a bilipschitz homeomorphism with respect to the quasihyperbolic metrics, see \cite{Vai-5}. Each conformal mapping between planar proper domains are quasihyperbolic maps. Moreover, an $L$-quasihyperbolic map is freely $\varphi$-quasiconformal with $\varphi(t)=L't$ for some $L'=L'(L)\geq 1$.

With the aid of Theorems \ref{main-thm-1} and \ref{main-thm-2}, we answer Question \ref{q-3} under some extra assumptions. One of our answers to Question \ref{q-3} is as follows.

\begin{theorem}\label{main-thm-3} Let $c,C,L\geq 1$ and  $D\subset E$, $D'\subset E'$ be domains. Suppose that
\begin{enumerate}
\item $f:\,\overline{D}\to \overline{D'}$ is a homeomorphism,
\item $f:\,D\to D'$ is $L$-quasihyperbolic, and $f:\,\partial D\to \partial D'$ is $L$-bilipschitz,
\item $D$ is $c$-uniform and $\partial D$ is $C$-uniformly perfect.
\end{enumerate}
Then $f:\,\overline{D}\to \overline{D'}$ is $L'$-bilipschitz with $L'$ depending only on $c,C$ and $L$.
\end{theorem}

Note that Theorem \ref{main-thm-3} is a partial answer to Question \ref{q-3} in requiring the domain $D$ to be $c$-uniform. However, the connectedness of $\partial D$ is replaced by a (weaker) uniformly perfect condition. Also, we do not require the domain $D$ to be bounded. On the other hand, this result is an improvement of \cite[Theorem 1.4]{LVX} with the absence of the uniformity of the target domain $D'$.

Furthermore, another answer to Question \ref{q-3} is Theorem \ref{main-thm-4} below.
Compared to Theorem \ref{main-thm-3}, the uniformity condition for domain $D$ is replaced by a relative quasim\"obius assumption for the mapping $f$.

\begin{theorem}\label{main-thm-4} Let $L,C\geq 1$,  $\theta:\,[0,\infty)\to[0,\infty)$ be a homeomorphism, and $D\subset E$, $D'\subset E'$ be domains. Suppose that $f:\, \overline{D}\to \overline{D'}$ is a homeomorphism, and that
\begin{enumerate}
\item $f:\,D\to D'$ is $L$-quasihyperbolic and $f:\,\partial D\to \partial D'$ is $L$-bilipschitz,
\item $f:\,\overline{D}\to \overline{D'}$ is $\theta$-quasim\"{o}bius relative to $\partial D$,
\item $\partial D$ is $C$-uniformly perfect.
\end{enumerate}
Then $f:\,\overline{D}\to \overline{D'}$ is $L'$-bilipschitz with $L'$ depending only on $C,L$  and $\theta$.
\end{theorem}

Finally, we study the relationship between quasihyperbolic maps and $(L,\vartheta)$-locally bilipschitz maps (See Definition \ref{de2-6}) which were introduced by Bonk et al. \cite[(1.20)]{BHK}. This class of mappings plays an important role in their uniforming procedure of Gromov hyperbolic spaces. Accordingly, it was shown in \cite{BHK} that there is a one-to-one correspondence between the bilipschitz classes of proper geodesic roughly starlike Gromov hyperbolic spaces and the quasisimilarlity classes of bounded locally compact uniform spaces, where quasisimilarity means quasisymmetry and $(L,\vartheta)$-locally bilipschitz. Recently, the first author in \cite{Z} established several equivalent conditions for quasihyperbolic maps associated to the $(L,\vartheta)$-locally bilipschitz condition. Using the results in \cite{Z}, we obtain the following result related to Question \ref{q-3}.

\begin{theorem}\label{main-thm-5}
Let $L,C\geq 1$, $\vartheta\in(0,1)$, and   $\theta:\,[0,\infty)\to[0,\infty)$ be a homeomorphism, and   $D\subset E$, $D'\subset E'$ be domains. Suppose that $f:\,\overline{D}\to \overline{D'}$ is a homeomorphism, and that
\begin{enumerate}
\item\label{z-00} $f:\,D\to D'$ is $\theta$-quasim\"{o}bius and $(L,\vartheta)$-locally bilipschitz,
\item $f:\,\partial D\to \partial D'$ is $L$-bilipschitz,
\item $\partial D$ is $C$-uniformly perfect.
\end{enumerate}
Then $f:\,\overline{D}\to \overline{D'}$ is $L'$-bilipschitz with $L'$ depending only on $C,L,\vartheta$ and $\theta$.
\end{theorem}

The paper is organized as follows. In Section \ref{sec-2}, we recall some basics about quasihyperbolic geometry, Gromov hyperbolic domains and mappings on Banach spaces. Section \ref{sec-3} is devoted to the proofs of Theorems \ref{main-thm-1} and \ref{main-thm-2}. The issue of Section \ref{sec-4} is to  answer Question \ref{q-3}. In particular, we prove Theorems \ref{main-thm-3}, \ref{main-thm-4} and \ref{main-thm-5}.

\section{Preliminaries and discussions}\label{sec-2}

\subsection{\bf Notations.} The basic notations are fairly standard and the same as in \cite{Vai-5}. We let $(X, |\cdot|)$ be a metric space, and
\[
B(x,r)=\{ z\in X:\; |z-x|<r\},\; \overline{B}(x,r)=\{ z\in X:\; |z-x|\leq r\},
\]
and $S(x,r)=\{ z\in X:\; |z-x|=r\}$. For a bounded set $A$ in $X$, $\diam (A)$ is the diameter of $A$.

For $C\geq 1$, a metric space $X$ is said to be {\it $C$-uniformly perfect}, if for each $x\in X$ and every $r>0$, $B(x,r)\setminus B(x, r/C)\not=\emptyset$ provided $X\setminus B(x,r)\not=\emptyset$. Note that each connected metric space is uniformly perfect.

For real numbers $s,t$, we set
\[ 
s\vee t=\max\{s,t\}\;\;\mbox{and}\;\;s\wedge t=\min\{s,t\}.
\]

To simplify expressions we often omit parentheses by writing $f(x)=x'$ etc.

\subsection{\bf Quasihyperbolic geometry and uniform domains.}
Throughout this paper, we use $E$ and $E'$ to denote real
Banach spaces with dimension at least $2$. The norm of a vector $x$
in $E$ is written as $|x|$, and for every pair of points $z_1$,
$z_2$ in $E$, the distance between them is denoted by $|z_1-z_2|$,
the closed line segment with endpoints $z_1$ and $z_2$ is denoted by $[z_1,
z_2]$.

Moreover, the one-point extension of $E$ is the Hausdorff space $\dot{E}=E\cup\{\infty\}$, where the neighborhoods of
$\infty$ are the complements of closed bounded sets of $E$. The boundary $\partial A$ and the closure $\overline{A}$ of a set $A\subset E$ are taken in $\dot{E}$. A proper
sub domain $D \subset E$ is a nonempty connected open set with $\partial D\neq \emptyset$, $d_D(x):=\dist(x,\partial D)$ for $x\in D$.

The {\it quasihyperbolic length} of a rectifiable arc or a path
$\alpha$ in the norm metric in a domain $D\subset E$ is the number (cf.
\cite{GO,GP}):
\[ 
\ell_k(\alpha)=\int_{\alpha}\frac{|dz|}{d_{D}(z)}.
\]

For each pair of points $x$, $y$ in $D$, the {\it quasihyperbolic  metric}
$k_D(x,y)$ between $x$ and $y$ is defined by
\[
k_D(x,y)=\inf\ell_k(\alpha),
\]
where the infimum is taken over all rectifiable curves $\alpha$ in $D$
joining $x$ to $y$.

For each pair of points $x$, $y$ in $D$, the {\it distance ratio  metric} $j_D(x,y)$ between $x$ and $y$ is defined by
\[ 
j_D(x,y)=\log\Big(1+\frac{|x-y|}{d_D(x)\wedge d_D(y)}\Big).
\]

For all $x$, $y$ in $D$, we have the following elementary inequalities (cf. \cite{Vai-5})
\[ 
k_{D}(x, y)\geq j_D(x, y)\geq \Big|\log \frac{d_{D}(x)}{d_{D}(y)}\Big|.
\]

Let $c\geq1$. A domain $D$ in $E$ is called $c$-{\it
uniform}, if each pair of points $x,y$ in $D$ can
be joined by a rectifiable arc $\alpha$ in $ D$ satisfying
\begin{enumerate}
\item $\ell (\alpha [x, z])\wedge \ell(\alpha[z, y])\leq c\,d_{D}(z)
$ for all $z\in \alpha$, and

\item $\ell(\alpha)\leq c\,|x-y|$,
\end{enumerate}
\noindent where $\ell(\alpha)$ denotes the length of $\alpha$ and
$\alpha[x,z]$ is the part of $\alpha$ between $x$ and $z$. Moreover, $\alpha$ is said to be a {\it uniform arc}.

We conclude this subsection with the following lemma for later use.

\begin{lemma}\label{zz1}$($\cite[Theorem 3.9]{Vai-5}$)$
Suppose that $D\subset E$ is a domain, and that $x$, $y\in D$ and either $|x-y|\leq \frac{1}{2}d_D(x)$ or $k_D(x,y)\leq 1$. Then
\[
\frac{1}{2}\frac{|x-y|}{d_D(x)}< k_D (x,y) \leq 2\frac{|x-y|}{d_D(x)}.
\]
\end{lemma}

\subsection{\bf Gromov hyperbolic domains.}
Assume that $D  \subset E$ is a domain and $\delta\geq 0$ is a constant. We say that $D$ is a {\it Gromov hyperbolic} or {\it $\delta$-hyperbolic} domain if
$(D ,k_{D})$ is a $\delta$-hyperbolic space, where $k_D$ is the quasihyperbolic metric of $D$. This means that
\[ 
(x|y)_p \geq (x|z)_p\wedge (z|y)_p-\delta
\] 
for all $x,$ $y,$ $z,$ $p\in D$, where $(x|y)_p$ is the Gromov product defined by
\[ 
2(x|y)_p =k_{D}(x,p)+k_{D}(y,p)-k_{D}(x,y).
\]

We recall some basic notations for a hyperbolic domain $D$. See \cite{Vai,Vai-00} for more details.
Fix a base point $p\in D$.
\begin{definition}
\begin{enumerate}
\item[{\rm (1)}] A sequence $\{x_i\}$ in $D$ is called a {\it Gromov sequence} if $(x_i|x_j)_p\rightarrow \infty$ as $i,$ $j\rightarrow \infty.$

\item[{\rm (2)}] Two such Gromov sequences $\{x_i\}$ and $\{y_j\}$ are said to be {\it equivalent} if $(x_i|y_i)_p\rightarrow \infty$ as $i\to \infty$.

\item[{\rm (3)}] The {\it Gromov boundary} $\partial^* D$ of $D$ is defined to be the set of all equivalence classes of Gromov sequences, and $D^*=D  \cup \partial^* D$ is called the {\it Gromov closure} of $D$.

\item[{\rm (4)}] For $\xi,$ $\zeta\in \partial^*D $, the Gromov product $(\xi|\zeta)_p$ of $\xi$ and $\zeta$ is defined by
\[ 
(\xi|\zeta)_p= \inf \big\{ \liminf_{i\rightarrow \infty}(x_i|y_i)_p:\; \{x_i\}\in \xi\;\;{\rm and}\;\; \{y_i\}\in \zeta\big\}.
\]
\end{enumerate}
\end{definition}

Assume that $D$ is a $\delta$-hyperbolic domain with $\delta>0$.
Let $0<\varepsilon\leq 1\wedge (5\delta)^{-1}$ and define
\[ 
\rho_{p,\varepsilon}(x,y)=e^{-\varepsilon(x|y)_p}
\]
for $x,y\in D^*$ with  the  convention $e^{-\infty}=0$. Because $\varepsilon \delta\leq 1$, we know from \cite[Proposition 5.16]{Vai-00} that there is a function $d_{p,\varepsilon}$ such that
\[
d_{p,\varepsilon}\leq \rho_{p,\varepsilon}\leq 2d_{p,\varepsilon},
\]
where the function $d_{p,\varepsilon}$ is a metametric on $D^*$, that is, it satisfies the axioms of a metric except that $d_{p,\varepsilon}(x,x)$ may be positive. In fact, $d_{p,\varepsilon}(x,y)=0$ if and only if $x=y\in \partial^* D$. Hence $d_{p,\varepsilon}$ defines a metric  on $\partial^* D$. We always say that $d_{p,\varepsilon}$ is a {\it visual metric} with base point $p\in D$ and parameter $\varepsilon>0$.

The metametric $d_{p,\varepsilon}$ defines a topology $\mathcal{T}^*$ in $D^*$. In this topology, the points of $D$ are isolated, and a sequence $\{x_i\}$ converges to a point $\xi\in \partial^* D$ if and only if $\{x_i\}$ is a Gromov sequence and $\{x_i\}\in \xi$. Since $\mathcal{T}^*$ in $D^*$ is discrete, the identity map $\id:\,D  \rightarrow D$ is continuous from the topology $\mathcal{T}^*$ to the metric topology. If it has a continuous extension $\varphi:\, D ^*\rightarrow \overline{D}$, we say that $\varphi$ is a {\it natural map}.

It was shown by Bonk et al. that a bounded domain in $\mathbb{R}^n$ is uniform if and only if it is both Gromov hyperbolic and its Euclidean boundary is naturally quasisymmetrically equivalent to the Gromov boundary, see \cite[Theorem 1.11]{BHK}. Subsequently, V\"ais\"al\"a obtained a dimension-free result which we recall now.

\begin{theorem}\label{Thm-D} $($\cite[Theorems 2.29 and 3.27]{Vai}$)$
Let $D$ be a proper domain of a Banach space $E$. Then $D$ is $c$-uniform if and only if it is $\delta$-hyperbolic and there is a natural $\theta$-quasim\"obius homeomorphism $\partial^*D\to \partial D$, where $c$ and $\delta$, $\theta$ depend only on each other.
\end{theorem}

\subsection{\bf  Quasisymmetric and quasim\"obius maps.}\label{sub-1}
 Let $X$ be a metric space and $\dot{X}=X\cup \{\infty\}$. By a
triple in $X$ we mean an ordered sequence $T=(x,y,z)$ of three
distinct points in $X$. The ratio of $T$ is the number
\[
\rho(T)=\frac{|y-x|}{|z-x|}.
\]
If $f:\, X\to Y$ is  an injective map, then the image of a triple  $T=(x,y,z)$  is the triple
$f(T)=(x',y',z')$.

Suppose that  $A\subset X$. A triple  $T=(x,y,z)$ in $X$ is said to
be a triple in the pair $(X, A)$ if either $x\in A$  or  $\{y,z\}\subset
A$. Equivalently, both $|y-x|$ and $|z-x|$ are distances from a
point in $A$.

\begin{definition} \label{def1-0} Let $X$ and $Y$ be metric spaces, and
$\eta:\, [0, \infty)\to [0, \infty)$ be a homeomorphism. Suppose that
$A\subset X$. An embedding $f:\, X\to Y$ is said to be {\it
$\eta$-quasisymmetric} relative to $A$, if $\rho(f(T))\leq \eta(\rho(T))$  for each triple $T$ in
$(X,A)$. \end{definition}

It is known that an embedding $f:\, X\to Y$ is $\eta$-quasisymmetric relative to $A$ if
and only if $\rho(T)\leq t$ implies that $\rho(f(T))\leq \eta(t)$
for each triple $T$ in $(X,A)$ and $t\geq 0$ (cf. \cite{Vai-5}).
Obviously, ``quasisymmetric relative to $X$" is equivalent to ordinary
``quasisymmetric".

If $\eta(t)=M(t^{1/\alpha}\vee t^{\alpha})$ for some constants $\alpha \geq 1$ and $M>0$, we call $f$ a {\it power quasisymmetric map}. We also record the following useful result.

\begin{theorem}\label{Thm-F} $($\cite[Theorems 3.5]{TV}$)$
Let $\eta:\, [0, \infty)\to [0, \infty)$ be a homeomorphism, and let $f:\,X\to Y$ be an $\eta$-quasisymmetric homeomorphism between bounded metric spaces. Then for $x\neq y\in X$ we have
\[ 
\frac{\diam(Y)}{2\eta(\frac{\diam (X)}{|x-y|})}\leq |x'-y'|\leq \eta\Big(\frac{2|x-y|}{\diam (X)}\Big)\diam (Y).
\]
\end{theorem}

A quadruple in $X$ is an ordered sequence $Q=(x,y,z,w)$ of four
distinct points in $X$. The cross ratio of $Q$ is defined to be the
number
\[
\tau(Q)=\tau(x,y,z,w)=\frac{|x-z|}{|x-y|}\cdot\frac{|y-w|}{|z-w|}.
\]
 Observe that this definition is extended  in usual sense if one of the points is
$\infty$. For example,
\[
\tau(x,y,z,\infty)= \frac{|x-z|}{|x-y|}.
\]
If $X_0 \subset \dot{X}$ and  $f:\, X_0\to \dot{Y}$
is an injective map,  then the image of a quadruple $Q$ in $X_0$ is the
quadruple $f(Q)=(x',y',z',w')$. Suppose that $A\subset X_0$. We say
that a quadruple $Q=(x,y,z,w)$ in $X_0$ is a quadruple in the pair
$(X_0, A)$ if  either $\{x,w\}\subset A$ or $\{y,z\}\subset A$.
Equivalently, all four distances in the definition of $\tau(Q)$ are
(at least formally) distances from a point in $A$.

\begin{definition} \label{def2'}
Let $X$ and $Y$ be metric spaces and   $\eta:\, [0, \infty)\to [0, \infty)$ be a
homeomorphism. Suppose that $A\subset \dot{X}$. An embedding $f:\,
\dot{X}\to \dot{Y}$ is said to be {\it $\eta$-quasim\"obius}
relative to $A$, if the inequality
$\tau(f(Q))\leq \eta(\tau(Q))$ holds for each quadruple in $(X,A)$.
\end{definition}

Apparently, ``$\eta$-$QM$ relative to $X$" is equivalent to ordinary ``quasim\"obius".

We conclude this subsection with the following useful notation, which was introduced by Bonk and  Kleiner in \cite{BK1}:
\[
\langle x,y,z,w \rangle=\frac{|x-z|\wedge |y-w|}{|x-y|\wedge |z-w|}.
\]
They also established a relation between $\tau(x,y,z,w)$ and $\langle x,y,z,w \rangle$ in the setting of metric spaces (\cite[Lemma 3.3]{BK1}).

\begin{lemma}\label{zz-1}
For any $x,$ $y,$ $z,$ $w$ in a metric space $X$, we have
\[
\langle x,y,z,w \rangle \leq \theta_0(\tau(x,y,z,w)),
\]
where $\theta_0(t)=3(t\vee \sqrt{t})$.
\end{lemma}

\subsection{\bf Mappings on domains.}\label{sub-3}
Let $C\geq 0$ and $M\geq 1$.
We say that a homeomorphism $f :\,D\to D'$ is {\it $C$-coarsely $M$-quasihyperbolic}, or briefly $(M, C)$-CQH, if it
satisfies
\[
M^{-1}k_D(x, y)-C\leq k_{D'}(x', y')\leq M k_D(x, y)+C
\]
for all $x$, $y\in D$.

\begin{definition}
We say that a homeomorphism $f:\, D\to D'$ is {\it $M$-bilipschitz} if
\[
M^{-1}\,|x-y| \leq |x'-y'|\leq M\,|x-y|
\]
for all $x$, $y\in D$,  and {\it $M$-QH} if
\[ 
M^{-1}\,k_{D}(x,y)\leq k_{D'}(x',y')\leq M\,k_{D}(x,y)
\]
for all $x$, $y\in D$.
\end{definition}

We remark that, if $f$ is $M$-bilipschitz or $M$-QH, then $f^{-1}$
has the same property. Now, we review the definitions of semi-solid, solid and FQC mappings introduced by V\"ais\"al\"a \cite{Vai-1,Vai-5}.

\begin{definition}\label{def-fqc}
Let $\varphi:\,[0,\infty)\to [0,\infty)$ be a homeomorphism with $\varphi(t)\geq t$. We say
that a homeomorphism $f:\, D\to D'$ is {\it $\varphi$-semisolid} if
\[
k_{D'}(x',y')\leq \varphi(k_{D}(x,y))
\]
for all $x$, $y\in D$, and {\it $\varphi$-solid} if both $f$ and $f^{-1}$ satisfy this condition.

We say that $f$ is {\it fully $\varphi$-semisolid}
(resp. {\it fully $\varphi$-solid}) if $f$ is
$\varphi$-semisolid (resp. $\varphi$-solid) on every  subdomain of $D$. In particular,
when $D=E$, the corresponding subdomains are taken to be proper ones. Fully $\varphi$-solid maps are also called {\it freely
$\varphi$-quasiconformal maps}, or briefly {\it $\varphi$-FQC maps}.
\end{definition}

Finally, we recall the definition of $(L,\vartheta)$-locally bilipschitz maps, which was introduced by Bonk et al. in \cite[(1.20)]{BHK}.

\begin{definition}\label{de2-6}
A homeomorphism $f:\,D\to D'$ is said to be {\it $(L,\vartheta)$-locally bilipschitz}, where $L\geq 1$ and $\vartheta\in(0,1)$, if for each $x\in D$ there is
a $c_x>0$ $($depending on $x)$ such that
\[
\frac{c_x}{L}|y-z|\leq |y'-z'|\leq Lc_x|y-z|,
\]
whenever $y,z\in B(x,\vartheta d_D(x))$.
\end{definition}

\section{Proofs of Theorems \ref{main-thm-1} and \ref{main-thm-2} and some applications}\label{sec-3}

\subsection{}
The purpose of this part is to show Theorem \ref{main-thm-1}. It is not hard to see that the sufficiency part follows from Theorem \ref{Thm-B}. Originally, to get Theorem \ref{Thm-B}, V\"ais\"al\"a made effective use of many tools and results in his study of the free quasiworld: the existence of neargeodesics, a notion of coarse length, the CQH invariance of solid arcs, and the length cigar and diameter cigar theorems for uniform domains. See \cite{Vai-2} for related definitions and discussions.

 As a reviewer of \cite{Vai-5}, Heinonen asked how to relate the freely quasiconformal mapping theory to the theory of Gromov hyperbolic spaces, see also \cite{Bu}. In view of this consideration, we provide a new proof  of Theorem \ref{Thm-B}. Our method is based on the Gromov hyperbolic characterization of uniform domains, see Theorem \ref{Thm-D}.

\subsection*{Proof of Theorem \ref{main-thm-1}.}
Suppose that $D$ is a $c$-uniform domain and $f:\,D\to D'$ is $(M, C)$-CQH.
We first prove the sufficiency  part and assume that $D'$ is $c'$-uniform. Since $D$ and $D'$ are uniform, it follows from \cite[Theorem 2.12]{Vai} that $(D,k_D)$ and $(D',k_{D'})$ are length $\delta$-hyperbolic spaces with $\delta\geq 0$ depending only on $c$ and $c'$, where $k_D$ and $k_{D'}$ are the quasihyperbolic metrics of $D$ and $D'$, respectively. Since $f$ is $(M,C)$-CQH, $f$ is a roughly quasi-isometric homeomorphism in the quasihyperbolic metrics. Thus we see from \cite[Theorem 5.38]{Vai-00} that $f$ has a
 homeomorphic extension to the Gromov closure
\[
f^*:\,D^*\to D'^*
\]
such that $f^*|_{\partial^*D}$ is $\theta^*$-quasim\"{o}bius with respect to any visual metrics $d_{p,\varepsilon}$ and $d_{q',\varepsilon'}$ for some positive parameters $\varepsilon=\varepsilon(\delta)$ and $\varepsilon'=\varepsilon'(\delta)$, where the control function $\theta^*$ depends only on $\delta,C,M$, but does not depend on the base points $p$ and $q'$.

On the other hand, since $D$ is $c$-uniform, with the aid of \cite[Proposition 2.26]{Vai}, we know that there is a natural bijective map
\[
\varphi:\,D^* \to \overline{D}
\]
such that every Gromov sequence $\{x_n\}\in \xi\in\partial^*D$ if and only if the sequence $\{x_n\}\subset D$ converges to a point $a=\varphi(\xi)\in\partial D$ in norm.

Moreover, by Theorem \ref{Thm-D}, it follows that the natural map
\[
\varphi:\,(\partial^*D,d_{p,\varepsilon}) \to (\partial{D},|\cdot|)
\]
is $\theta_1$-quasim\"{o}bius for some control function $\theta_1$ depending only on $c$ and $\delta$.

Similarly, since $D'$ is $c'$-uniform, there is a natural bijective map
\[
\psi:\,D'^* \to \overline{D'}
\]
such that every Gromov sequence $\{x_n'\}\in \xi'\in\partial^*D'$ if and only if the sequence $\{x_n'\}\subset D'$ converges to a point $a'=\psi(\xi')\in\partial D'$ in norm; and moreover,
\[
\psi:\,(\partial^*D',d_{q',\varepsilon'}) \to (\partial{D'},|\cdot|)
\]
is $\theta_2$-quasim\"{o}bius for some control function $\theta_2$ depending only on $c'$ and $\delta$.

So we see that
\[
\overline{f}:\,=\psi\circ f^* \circ \varphi^{-1}:\, \overline{D}\to \overline{D'}.
\]
is a homeomorphic extension of $f$ and   is also $\theta$-quasim\"{o}bius on $\partial D$ with
$\theta=\theta_2\circ \theta^* \circ \theta_1^{-1},$
because the composition and  the inverse of quasim\"{o}bius maps are also quasim\"{o}bius, quantitatively. This proves the sufficiency.

It remains to show the necessity. Assume that $f$ extends to a homeomorphism $\overline{f}:\,\overline{D}\to \overline{D'}$ and $\overline{f}$ is $\theta$-quasim\"obius on $\partial D$.
We need to show that $D'$ is $c'$-uniform. To this end, by Theorem \ref{Thm-D}, we only need to show that there is a constant $\delta'$ and a homeomorphism $\theta':\,[0.\infty)\to [0,\infty)$ depending only on the data of $D$ and $f$ such that
\begin{enumerate}
  \item\label{l-1} $(D',k_{D'})$ is $\delta'$-hyperbolic;
  \item\label{l-2} The identity map $\id_{D'}:\, D'\to D'$ has a natural extension $\phi:\,D'^* \to \overline{D'}$;
  \item\label{l-3} The natural map $\phi:\,(\partial^*D',d_{q',\varepsilon'}) \to (\partial{D'},|\cdot|)$ is $\theta'$-quasim\"obius, where $d_{q',\varepsilon'}$ is a visual metric with parameter $\varepsilon'=\varepsilon'(\delta')$ and a base point $q'\in D'$.
\end{enumerate}

We first verify (\ref{l-1}). From \cite[Theorem 2.12]{Vai} it follows that $(D,k_D)$ is a length $\delta$-hyperbolic space where $\delta$ depends only on $c$. Then we know from \cite[Theorem 3.20]{Vai-00} that $(D',k_{D'})$ is a length $\delta'$-hyperbolic space with $\delta'$ depending only on $c,C$ and $M$, because $f$ is assumed to be $(M,C)$-CQH. Hence we are done.

To prove (\ref{l-2}), we immediately see from the proof of the sufficiency part that
\begin{enumerate}
  \item[(a)] $f$ has a homeomorphic extension to the Gromov closure $f^*:\,D^*\to D'^*$ such that $f^*|_{\partial^*D}:=\partial f^*$ is $\theta^*$-quasim\"{o}bius with respect to any visual metrics $d_{p,\varepsilon}$ and $d_{q',\varepsilon'}$ for some positive parameters $\varepsilon=\varepsilon(\delta)$ and $\varepsilon'=\varepsilon'(\delta')$, where the control function $\theta^*$ depends only on $\delta,\delta',C$ and $M$, but does not depend on the base points $p\in D$ and $q'\in D'$.
  \item[(b)]  There is a natural homeomorphism $\varphi:\,D^* \to \overline{D}$ such that $\varphi:\,(\partial^*D,d_{p,\varepsilon}) \to (\partial{D},|\cdot|)$ is $\theta_1$-quasim\"{o}bius for some control function $\theta_1$ depending only on $c$ and $\delta$.
\end{enumerate}
Thus $f$ induces a mapping
\[
\phi:=\overline{f}\circ \varphi\circ (f^*)^{-1}:\, D'^*\to \overline{D'}.
\]
Since $f$ extends to a homeomorphism $\overline{f}:\,\overline{D}\to \overline{D'}$, we know from the facts (a) and (b) that $\phi$ is a homeomorphism. Moreover, since $\varphi|_D=\id_D$, we have $\phi|_{D'}=\id_{D'}$. Hence $\phi$ is a natural map, as desired.

It remains to show (\ref{l-3}). Since $\overline{f}$ is $\theta$-quasim\"obius on $\partial D$, it follows immediately from the facts (a) and (b) that
\[
\phi:\, (\partial^* D',d_{q',\varepsilon'})\to (\partial{D'},|\cdot|)
\]
is $\theta'$-quasim\"obius with $\theta'=\theta\circ \theta_1\circ (\theta^*)^{-1}$, because the composition and  the inverse of quasim\"{o}bius maps are also quasim\"{o}bius, quantitatively.

Hence the necessity follows and the proof of this theorem is complete.
\qed

\subsection{}
Our aim of this subsection is to prove Theorem \ref{main-thm-2}. Before the proof, we remark the sufficiency and the fact that quasisymmetry implies quasim\"obius are known, see \cite[Theorems 3.2 and 3.12]{Vai-0}. In this note, we also give different proofs of these facts by using a useful notion which was introduced by Bonk and Kleiner, see Subsection \ref{sub-1}. The original proof given by V\"ais\"al\"a  needed to consider  several cases. However, our proof indicates the precise dependence between the constant of the three-point condition, and the control functions of quasisymmetric and quasim\"obius mappings.

\subsection*{Proof of Theorem \ref{main-thm-2}.}
Suppose that $f:\,X\to Y$ is a homeomorphism between bounded connected metric spaces. We first show the sufficiency. Assume that $f$ is $\theta$-quasim\"{o}bius which satisfies the $\lambda$-three-point condition for some constant $\lambda\geq 1$ and points $z_1,$ $z_2,$ $z_3$ in $X$. Fix three distinct points $x,a,b\in X$.
It suffices to find some homeomorphism $\eta:\,[0,\infty)\to [0,\infty)$ such that
\be\label{l-4}
|x'-a'|\leq \eta(t)|x'-b'| ~\mbox{ whenever}\;\;|x-a|=t|x-b|.
\ee
First, by the $\lambda$-three-point condition, we see that for each $u\in X$, there are $i, j\in\{1,2,3\}$ with $i\neq j$ such that
\be\label{l-5}
|z_i-u|\wedge |z_j-u|\geq \frac{\diam (X)}{3\lambda}.
\ee

Similarly, there are $m, n\in\{1,2,3\}$ with $m\neq n$ such that for any $v'\in Y$,
\be\label{l-6}
|z_m'-v'|\wedge |z_n'-v'|\geq \frac{\diam (Y)}{3\lambda}.
\ee
Therefore, by (\ref{l-5}) and (\ref{l-6}), we obtain that there must exist a point $z_i\in \{z_1,z_2,z_3\}$ such that
\[
|z_i-a|\geq \frac{\diam (X)}{3\lambda} ~\mbox{ and }~ |z_i'-b'|\geq \frac{\diam (Y)}{3\lambda}
\]
and so
\[
|z_i-a|\wedge |x-b|\geq \frac{|x-b|}{3\lambda}\;\;\mbox{and}\;\; |x'-a'|\wedge |z_i'-b'|\geq \frac{|x'-a'|}{3\lambda}.
\]
Thus, we obtain
\be\label{l-7}
\langle x',b',a',z_i'\rangle\geq \frac{|x'-a'|}{3\lambda |x'-b'|}\;\;\mbox{and}\;\; \langle x,b,a,z_i\rangle\leq \frac{3\lambda |x-a|}{|x-b|}.
\ee

On the other hand, since $f$ is $\theta$-quasim\"{o}bius, we see from Lemma \ref{zz-1} that
\be\label{l-8}
\langle x',b',a',z_i'\rangle \leq \theta'(\langle x,b,a,z_i\rangle),
\ee
where $\theta'(t)=\theta_0\circ \theta\circ \Big(\frac{1}
{\theta_0^{-1}(1/ t)}\Big)$ and $\theta_0(t)=3(t\vee\sqrt{t})$.
Then we obtain by (\ref{l-7}) and (\ref{l-8}) that
\[
\frac{|x'-a'|}{|x'-b'|}\leq 3\lambda \theta'(\langle x,b,a,z_i\rangle)\leq 3\lambda \theta'\Big(3\lambda\frac{|x-a|}{|x-b|}\Big)=3\lambda \theta'(3\lambda t).
\]
Hence we get (\ref{l-4}) by taking $\eta(t)=3\lambda \theta'(3\lambda t).$ This ensures the sufficiency.

Next, we prove the necessity by assuming that $f$ is $\eta$-quasisymmetric. We first verify that $f$ is $\theta$-quasim\"obius. For any distinct points $x,y,z,w\in X$, no generality is lost in assuming that $|x'-y'|\leq |z'-w'|$. Since $f$ is $\eta$-quasisymmetric, we compute
\begin{eqnarray*} \langle x',y',z',w'\rangle
&=& \frac{|x'-z'|\wedge |y'-w'|}{|x'-y'|}\\
&\leq& \eta\Big(\frac{|x-z|}{|x-y|}\Big) \wedge \eta\Big(\frac{|y-w|}{|x-y|}\Big)
\\&=& \eta\Big(\frac{|x-z|\wedge |y-w|}{|x-y|}\Big)
\\&\leq& \eta(\langle x,y,z,w\rangle).
\end{eqnarray*}
Therefore, Lemma \ref{zz-1} ensures that $f$ is $\theta$-quasim\"obius with $\theta(t)=\theta_0^{-1}(\frac{1}{\eta\circ \theta_0(t)})^{-1}$.

It remains to verify the three-point condition for $f$. To this end, choose two points $z_1,z_2\in X$  such that $|z_1-z_2|\geq \diam (X)/2$. Then we claim that there is
a $z_3\in X$ such that for all $i=1,2$, we have
\be\label{z-04}
|z_3'-z_i'|\geq \diam (Y)/6.
\ee
Otherwise,
\[ 
Y=B_1\cup B_2,
\]
where $B_i=\overline{B}(z_i',\diam (Y)/6)$, for $i=1,2.$ First we observe that
\[ 
B_1\cap B_2\neq \emptyset,
\]
since $Y$ is connected. Thus we have
\[
|w_1'-w_2'|\leq |w_1'-z_1'|+|z_1'-z_2'|+|z_2'-w_2'|\leq \frac{2\diam (Y)}{3},
\]
for all $w_i'\in B_i$, $i=1,2$. From which it follows that
\[ 
Y\setminus (B_1\cup B_2)\neq \emptyset,
\]
which is absurd and so (\ref{z-04}) follows.

Finally, we shall find some constant $\lambda\geq 1$ such that $f$ satisfies the $\lambda$-three-point condition with respect to the points $z_1,$ $z_2,$ $z_3$. With the aid of Theorem \ref{Thm-F}, we find that
\[
|z_1'-z_2'|\geq \frac{\diam (Y)}{2\eta(\diam (X)/|z_1-z_2|)}\geq  \frac{\diam (Y)}{2\eta(2)}.
\]
Now by (\ref{z-04}) we have for all $i=1,2$,
\[ 
\frac{1}{6} \leq  \frac{|z_i'-z_3'|}{\diam (Y)}\leq \eta\Big(\frac{2|z_i-z_3|}{\diam (X)}\Big),
\]
which implies that
\[ 
|z_i-z_3|\geq \frac{1}{2}\eta^{-1}(1/6)\diam (X).
\]
The proof of the theorem is complete if we let $\lambda=\max\{6,2\eta(2),2/\eta^{-1}(1/6)\}.$
\qed

\subsection{Applications.}\label{sub-2}
In this subsection, as applications of Theorems \ref{main-thm-1} and \ref{main-thm-2}, we investigate the boundary behavior of freely quasiconformal mappings on
uniform domains of  Banach spaces. We begin with an elementary fact, Lemma \ref{z-1} below. We failed to find this fact explicitly  in the literature, although its content should be well-known; for the completeness we provide the details.

\begin{lemma}\label{z-1}  Let $D\subset E$ be a bounded domain. Then
\be\label{z-2}
\diam (\partial D)=\diam(D) = \diam (\overline{D})<\infty,
\ee
and for any $z_1,z_2\in\partial D$, there is a point $z_3\in\partial D$ such that
\be\label{z-33}
|z_1-z_2|\wedge |z_3-z_2|\geq \diam (D)/6.
\ee
\end{lemma}

\bpf Since $\diam (\partial D) \leq \diam (D) \leq \diam (\overline{D})<\infty$, to get $(\ref{z-2})$ we only need to show $\diam (\overline{D}) \leq \diam (\partial D)$. Let $\varepsilon>0$ be given small enough. Choose $x,y\in \overline{D}$ with $\diam (\overline{D})\leq |x-y|+\varepsilon$. If $x$ and $y$ are both contained in $\partial D$,
 then there is nothing to prove. Without loss of generality, we may assume that $x\in D$. Denote  by $L_{\overrightarrow{yx}}$, the ray emanating from $y$ through $x$ to $\infty$. Thus
\[
L_{\overrightarrow{yx}} \cap \partial D\neq \emptyset,
\]
since $D$ is bounded. Take $x_0\in L_{\overrightarrow{yx}} \cap \partial D$. If $y\in D$ and continue in this manner, then
 we find that there is a point $y_0\in L_{\overrightarrow{xy}} \cap \partial D$. For the other possibility that $y\in\partial D$, let $y_0=y$. Then we have
 \[
\begin{aligned}[t]
\diam (\overline{D}) \leq |x-y|+\varepsilon & 
= |x_0-y_0|-|x-x_0|-|y_0-y|+\varepsilon\\
& \leq |x_0-y_0|+\varepsilon\leq \diam (\partial D)+\varepsilon.
\end{aligned}
\]
Letting $\varepsilon\to 0$, $(\ref{z-2})$ follows.

It remains to verify $(\ref{z-33})$. To this end, assuming the opposite, and using $(\ref{z-2})$ we find that there are two points $z_1,z_2\in\partial D$ such that
\[
\partial D\subset B(z_1,\diam (D)/6) \cup B(z_2,\diam (D)/6) =:S.
\]
Since $D$ is bounded, we know that $\overline{D} \subset S$. We claim that
\be\label{z-01}
B(z_1,\diam (D)/6) \cap B(z_2,\diam (D)/6) \neq \emptyset.
\ee

Otherwise, let $x_n,y_n\in D$ be such that $x_n\to z_1$ and $y_n\to z_2$ as $n\to \infty$. For each $n$, choose a curve $\gamma_n$ joining $x_n$ and $y_n$ in $D$. Then there is some point
\[
w_n\in \gamma_n\cap (E\setminus S)\subset D\subset S,
\]
which is absurd and so (\ref{z-01}) is true.

Now from (\ref{z-01}) it follows that $|z_1-z_2|\leq \diam (D)/3$. Furthermore, for all $u_1\in B(z_1,\diam (D)/6)$ and $u_2\in B(z_2,\diam (D)/6)$, we have
\[
|u_1-u_2|\leq |u_1-z_1|+|z_1-z_2|+|z_2-u_2|\leq 2\diam (D)/3.
\]
Combining this estimate to the opposite assumption we obtain
\[
\diam (\partial D)\leq \diam (S)\leq 2\diam (D)/3,
\]
which clearly contradicts $(\ref{z-2})$. The proof of the lemma is complete.
\epf

Using Lemma \ref{z-1} and Theorem \ref{main-thm-2} we obtain the following result associated with the quasisymmetry of quasim\"{o}bius maps between domains in Banach spaces whenever their boundaries are quasisymmetrically equivalent.

\begin{theorem}\label{s-7}
Let $D\subset E$ and $D'\subset E'$ be domains, and  $f:\,\overline{D}\to \overline{D'}$ be a homeomorphism. If $f:\,D\to D'$ is  $\theta$-quasim\"{o}bius and $f:\,\partial D\to \partial D'$ is $\theta$-quasisymmetric, then $f:\,\overline{D}\to \overline{D'}$ is $\eta$-quasisymmtric with $\eta$ depending only on $\theta$.
\end{theorem}

\bpf
We divide the proof into two cases. Note that $\partial D$ and $\partial D'$ are taken in the Hausdorff spaces $\dot{E}$ and $\dot{E'}$, respectively. If $D$ is unbounded, then
\cite[Theorem 3.10]{Vai-0} ensures that $f(\infty)=\infty$, because $f|_{\partial D}$ is quasisymmetric. Again by \cite[Theorem 3.10]{Vai-0} and the assumption that $f:\,D\to D'$ is $\theta$-quasim\"{o}bius, we deduce that $f$ is $\theta$-quasisymmtric.

 We next consider the remaining case that $D$ is bounded. To prove that $f:\,\overline{D}\to \overline{D'}$ is $\eta$-quasisymmtric, by Theorem \ref{main-thm-2}, we only need to show the following:
\bcl  There are three points $z_i\in \overline{D}$, $i=1,2,3$, such that
\[
|z_i-z_j|\geq \frac{1}{\lambda}\diam (\overline{D})\;\;\mbox{and}\;\; |z'_i-z'_j|\geq \frac{1}{\lambda}\diam (\overline{D'})
\]
for all $i, j\in\{1,2,3\}$ with $i\neq j$ and some suitable constant $\lambda\geq 1$ decided later which depends only on $\theta$.
\ecl
To this end, Lemma \ref{z-1} ensures that there are two points $z_1,z_2\in \partial D$  such that
\be\label{ll-1}
|z_1-z_2|\geq \diam (\partial D)/2=\diam  (\overline{D})/2.
\ee
 Now using (\ref{z-33}), we obtain that there is a point $z_3\in \partial D$ such that
\be\label{ll-2}
|z_1'-z_2'|\wedge |z_3'-z_2'|\geq \diam (D')/6.
\ee

Since $f:\,\partial D\to \partial D'$ is $\theta$-quasisymmetric, by Theorem \ref{Thm-F}, we compute
\be\label{ll-3}
|z_1'-z_2'|\geq \frac{\diam (D')}{2\theta(\diam (D)/|z_1-z_2|)}\geq  \frac{\diam (D')}{2\theta(2)},
\ee
and
\[
1/6\leq  \frac{|z_i'-z_3'|}{\diam (D')}\leq \theta\Big(\frac{2|z_i-z_3|}{\diam (D)}\Big),\;\;\mbox{for\,\,\,all}\;\;\;i=1,2,
\]
which implies
\be\label{ll-4}
|z_i-z_3|\geq \frac{1}{2}\theta^{-1}(1/6)\diam (D).
\ee
Then by (\ref{ll-1}), (\ref{ll-2}), (\ref{ll-3}) and (\ref{ll-4}), the claim follows by choosing
\[
\lambda=\max\{6,2\theta(2),2/\theta^{-1}(1/6)\}.
\]
The proof of the theorem is complete.
\epf

From a similar argument as in the  proofs of Theorem \ref{s-7} and \cite[Theorem 6.29]{Vai-5}, it is not hard to get the following corollary.

\begin{corollary}\label{s-8}
Let $D\subset E$ and $D'\subset E'$ be domains. If $f:\,D\to D'$ is $\theta$-quasim\"{o}bius relative to $\partial D$ and $f:\,\partial D\to \partial D'$ is $\theta$-quasisymmetric, then $f:\,\overline{D}\to \overline{D'}$ is $\eta$-quasisymmtric relative to $\partial D$ with $\eta$ depending only on $\theta$.
\end{corollary}

With the aid of Theorems \ref{main-thm-1} and \ref{s-7}, we prove

\begin{theorem}\label{s-9}
Let $D\subset E$ and $D'\subset E'$ be domains, and $f:\,\overline{D}\to \overline{D'}$ be a homeomorphism. Suppose that $f:\,\partial D\to \partial D'$ is $\theta$-quasisymmetric  (resp. $\theta$-quasim\"{o}bius). Then $f:\,D\to D'$ is $\varphi$-FQC and $D$ is $c$-uniform if and only if $D'$ is $c'$-uniform and $f:\,\overline{D}\to \overline{D'}$ is $\eta$-quasisymmtric (resp. $\eta$-quasim\"{o}bius), where $\varphi$, $c$ and $\eta$, $c'$ depend only on each other and $\theta$.
\end{theorem}
\bpf Because each $\eta$-quasim\"{o}bius mapping between proper domains is $\varphi$-FQC with $\varphi$ depending only on $\eta$ (cf. \cite[Theorem 7.17]{Vai-5}), the sufficiency follows immediately from Theorems \ref{Thm-A} and \ref{s-7}. Thus, we only need to show the necessity in assuming that $f$ is $\varphi$-FQC and $D$ is $c$-uniform.

We first assume that $f:\,\partial D\to \partial D'$ is $\theta$-quasim\"{o}bius, and we shall show that $f:\,\overline{D}\to \overline{D'}$ is $\eta$-quasim\"{o}bius. Since $f$ is $\varphi$-FQC, we see from \cite[Theorem 4.7]{Vai-5} that $f$ is $(M,C)$-CQH for some constants $M$ and $C$ depending only on $\varphi$. Thus using this fact and Theorem \ref{main-thm-1}, one immediately sees that $D'$ is $c'$-uniform with $c'$ depending only on $c,\varphi$ and $\theta$. Moreover, it follows from the assumption that $f$ is $\varphi$-FQC and Theorem \ref{Thm-A} that $f:\,D\to D'$ is $\eta$-quasim\"{o}bius with $\eta$ depending only on $c,\varphi$ and $\theta$. The assertion follows in this case by \cite[Theorem 3.18]{Vai-0}.

Next we assume that $f:\,\partial D\to \partial D'$ is $\theta$-quasisymmetric. We need to verify the quasisymmetry of $f:\,\overline{D}\to \overline{D'}$. Since every quasisymmetric mapping is quasim\"{o}bius, by the previous arguments, we see that $f:\,D\to D'$ is $\theta'$-quasim\"{o}bius with $\theta'$ depending only on $c,\varphi$ and $\theta$. With the help of this fact and Theorem \ref{s-7}, the quasisymmetry of $f$ follows. This completes the proof.
\epf

\begin{remark}
Theorem \ref{s-9} is an improvement of \cite[Theorem 1.2]{LVX} under weaker boundary hypotheses for $f$ and the absence of the uniformity of the target domain $D'$.
\end{remark}

Finally, we consider the boundary behavior of  freely quasiconformal mappings on uniform domains of the one-point extension space $\dot{E}=E\cup\{\infty\}$. Note that in the case that a domain $D\subset \dot{E}$, if $D$ is bounded, then $D$ is also a domain of $E$. If $D$ is unbounded, then we may have
 either $\infty\in D$ or $\infty\in \partial D$. This is the only difference with the case that a domain $D\subset E$.

\begin{theorem}\label{s-10}
Let $D\subset \dot{E}$ and $D'\subset \dot{E'}$ be domains, and  $f:\,\overline{D}\to \overline{D'}$ be a homeomorphism with $f:\,\partial D\to \partial D'$ being $\theta$-quasisymmetric. Then $f:\,D\to D'$ is  $\varphi$-FQC which maps bounded sets onto bounded sets, and $D$ is $c$-uniform if and only if $D'$ is $c'$-uniform and $f:\,\overline{D}\to \overline{D'}$ is $\eta$-quasisymmetric, quantitatively.
\end{theorem}
\bpf
It is known that quasisymmetric mappings map bounded sets onto bounded sets, see \cite{TV}. Moreover, one observes from \cite[Theorem 7.12]{Vai-5} that $\eta$-quasisymmetric mappings between proper domains are $\varphi$-FQC with $\varphi$ depending only on $\eta$. So the sufficiency follows from Theorem \ref{s-9}.

We only need to show the necessity. Assume that $f:\,D\to D'$ is $\varphi$-FQC which maps bounded sets onto bounded sets and that $D$ is a $c$-uniform domain. We again see from Theorem \ref{s-9} that $D'$ is $c'$-uniform. It remains to show that $f$ is quasisymmetric. The proof of this assertion is divided naturally into three parts.

If $D$ is bounded, then clearly $D\subset E$ is a domain. Using the hypothesis that $f$ maps bounded sets onto bounded sets, we obtain that $D'\subset E'$ is also a bounded domain.
The assertion follows directly from Theorem \ref{s-9}.

On the other hand, if $D$ is unbounded and if $\infty \not\in D$, then in this case $D$ again belongs to $E$ and $\infty\in\partial D$. Since $f$ is quasisymmetric on $\partial D$, we have $f(\infty)=\infty$. The same reason as the previous argument gives the assertion.

Finally,  suppose that $D$ is unbounded and $\infty \in D$. Since $f$ maps bounded sets onto bounded sets, we see that $f(\infty)=\infty$. Hence $f:\,\overline{D}\to \overline{D'}$ is $\theta$-quasisymmetric by means of Theorem \ref{s-9} and \cite[Theorem 3.10]{Vai-0}. Therefore, the theorem follows.
\epf

 We conclude this section with an example to demonstrate that the condition ``$f$ maps bounded sets onto bounded sets" of Theorem \ref{s-10} cannot be dropped.

\beg\label{exa-1}
Let $D=\{z\in \mathbb{C}:\, |z|<1\}$ and $D'=u(D)=\mathbb{C}\cup\{\infty\}\setminus \overline{D}$, where $u(z)=1/z$ for all $z\in \mathbb{C}\cup \{\infty\}$. Thus $u$ is a M\"obius transformation with identity boundary value but not quasisymmetric.
\eeg

\section{Boundary behavior of quasihyperbolic maps in Banach spaces}\label{sec-4}
\subsection{Quasihyperbolic maps in Banach spaces}
We begin this section with a natural question related to Theorem \ref{s-9}.

\bques\label{q-4}
Under the assumption of Theorem \ref{s-9}, suppose further that $f:\,\partial D\to \partial D'$ is $L$-bilipschitz.
 Is $f$ bilipschitz? In particular, is a quasim\"{o}bius map with bilipschitz boundary correspondence bilipschitz?
\eques

The following example tells us that the answer to this question is negative in general.

\beg\label{exa-2}
 Let $\alpha>1$, $E$ be a Banach space with $B=B(0,1)\subset E$,  and  $f:\,B\setminus \{0\}\to B\setminus \{0\}$ be the radial power map
\[
f(x)=|x|^{\alpha-1}x,\;\;\;\;\;\;x\in B\setminus \{0\}.
\]
It is known that $f$ is $\varphi_\alpha$-FQC and has identical boundary value, and is quasisymmetric, see \cite{Vai-5}. However, it is not bilipschitz.
\eeg

It turns out that some extra conditions are necessary to investigate bilipschitz boundary correspondence of freely quasiconformal mappings. Accordingly, V\"{a}is\"{a}l\"{a} raised Question \ref{q-3} for the class of quasihyperbolic maps. The aim of this section is to consider Question \ref{q-3}. In order to prove Theorems \ref{main-thm-3}, \ref{main-thm-4} and \ref{main-thm-5}, we need some more  preparation. We begin with the following equivalent conditions for quasihyperbolic maps.

\begin{lemma}$($\cite[Theorem 1.1]{Z}$)$\label{z-3}
Let $D\subset E$ and $D'\subset E'$ be domains, and  $f:\,D\to D'$ be a homeomorphism. Then the following conditions are quantitatively equivalent:
\begin{enumerate}
\item\label{qh-1.0} $f:\,D\to D'$ is $M$-QH;
\item\label{qh-1.2} Both $f$ and $f^{-1}$ are $(L,\vartheta)$-locally bilipschitz,
\item\label{qh-1.1} There are constants $C\geq 1$ and $\mu\in(0,1)$ such that
\[ 
\frac{1}{C}\frac{|x-y|}{d_D(x)}\leq \frac{|x'-y'|}{d_{D'}(x')}\leq C\frac{|x-y|}{d_D(x)},
\]
whenever $x,y\in D$ with $|x-y|<\mu d_D(x)$.
\end{enumerate}
\end{lemma}

The following auxiliary lemma asserts that the composition of a quasihyperbolic map and a locally bilipschitz map is locally bilipschitz, quantitatively.

\begin{lemma}\label{z-1.0}
Suppose that $E$, $T$ and $H$ are Banach spaces, each of which has dimension at least $2$, and that $D\subset E$, $G\subset T$, $\Omega\subset H$ are proper domains. Let $M,L\geq 1$ and $\vartheta\in(0,1)$, and let $f:\,D\to G$, $g:\,G\to \Omega$ be homeomorphisms. We have the following:
\begin{enumerate}
\item\label{z-1.1} If $f$ is $M$-QH and $g$ is $(L,\vartheta)$-locally bilipschitz, then $g\circ f$ is $(L',\vartheta')$-locally bilipschitz.
\item\label{z-1.2} If $f$ is $(L,\vartheta)$-locally bilipschitz and $g$ is $M$-QH, then $g\circ f$ is $(L'',\vartheta'')$-locally bilipschitz.
\end{enumerate}
Here $L'$, $L''$, $\vartheta'$, and $\vartheta''$ are constants depending only on $M,L$ and $\vartheta$.
\end{lemma}

\bpf Let $C\geq 1$ and $\mu\in(0,1)$ be constants given in Lemma \ref{z-3} which depends only on $M$.

We first prove (\ref{z-1.1}). Let $\vartheta'=\vartheta\mu/C$. Fix $x\in D$ and $y,z\in B(x,\vartheta' d_D(x))$.
Since $f$ is $M$-QH, Lemma \ref{z-3} guarantees that there are constants $L_0=L_0(M)$ and $c_x>0$ such that
\be\label{z-1.3}
\frac{c_x}{L_0}|y-z|\leq |f(y)-f(z)|\leq L_0c_x|y-z|.
\ee
On the other hand, again by Lemma \ref{z-3}, we find that
\[
\frac{|f(x)-f(u)|}{d_{G}(f(x))}\leq C \frac{|x-u|}{d_D(x)}<C\vartheta'\leq \vartheta ~\mbox{ for all $u\in B(x,\vartheta' d_D(x))$},
\]
which yields that $f(y),f(z)\in B(f(x),\vartheta d_G(f(x)))$. Moreover, since $g$ is $(L,\vartheta)$-locally bilipschitz, it follows that there is some $c_x'>0$ such that
\be\label{z-1.4}
\frac{c_x'}{L}|f(y)-f(z)|\leq |g\circ f(y)-g\circ f(z)|\leq L c_x'|f(y)-f(z)|.
\ee

Then by (\ref{z-1.3}) and (\ref{z-1.4}) we immediately find that $g\circ f$ is $(L',\vartheta')$-locally bilipschitz with the choice of $L'=LL_0$. This proves (1).

Next, we show (\ref{z-1.2}). Set $\vartheta''=\vartheta\mu/L^2$. We demonstrate that $g\circ f$ is $(L'',\vartheta'')$-locally bilipschitz with $L''=CL$.

Since $f$ is $(L,\vartheta)$-locally bilipschitz, a direct computation shows that $f$ is $\eta$-quasisymmetric on $B(x,\vartheta d_D(x))$ for all $x\in D$ with $\eta(t)=L^2t$. Moreover, by \cite[Theorem 6.12]{Vai-5} we may assume that $f$ is also $\eta$-quasisymmetric on $\overline{B}(x,\vartheta d_D(x))$. Then for all $v\in B(x,\vartheta d_D(x))$ with $|v-x|=\vartheta'' d_D(x)<\vartheta d_D(x)$ and for all $a\in S(x,\vartheta d_D(x))$, we have
\[
|f(x)-f(v)|\leq L^2 \frac{|x-v|}{|x-a|} |f(x)-f(a)|= \frac{L^2\vartheta''}{\vartheta} |f(x)-f(a)|.
\]
Now by the arbitrariness of $a\in S(x,\vartheta d_D(x))$ we get
\be\label{z-1.5}
|f(x)-f(v)|\leq \frac{L^2\vartheta''}{\vartheta} d_G(f(x)).
\ee

For each $x\in D$, since $f$ is $(L,\vartheta)$-locally bilipschitz, there is some $c_x>0$ such that
\be\label{z-1.6}
\frac{c_x}{L}|y-z|\leq |f(y)-f(z)|\leq Lc_x|y-z| ~\mbox{ for  all $y,z\in B(x,\vartheta' d_D(x))$.}
\ee

On the other hand, by (\ref{z-1.5}) and by the choice of $\vartheta''$, we see that $f(y),f(z)\in B(f(x),\mu d_G(f(x)))$. Moreover, since $g$ is $M$-QH, by Lemma \ref{z-3} it follows that
\be\label{z-1.7}
\begin{aligned}[t] 
|f(y)-f(z)|\frac{d_{\Omega}(g\circ f(x))}{C d_G(f(x))} & \leq |g\circ f(y)-g\circ f(z)|\\
&\leq |f(y)-f(z)|\frac{Cd_{\Omega}(g\circ f(x))}{d_G(f(x))}.
\end{aligned}
\ee

Finally, we see from (\ref{z-1.6}) and (\ref{z-1.7}) that $g\circ f$ is $(L'',\vartheta'')$-locally bilipschitz. This completes the proof.
\epf

We conclude this part by showing the following lemma which is needed in the proof of Theorem \ref{main-thm-5}.

\begin{lemma}\label{z-2.0}
Let $L\geq 1$, $\vartheta\in(0,1)$ and $\theta:\,[0,\infty)\to[0,\infty)$ be a homeomorphism, and $D\subset E$, $D'\subset E'$ be domains. If $f:\,D\to D'$ is $\theta$-quasim\"{o}bius and $(L,\vartheta)$-locally bilipschitz, then $f^{-1}$ is
$(L',\vartheta')$-locally bilipschitz with $L'$ and $\vartheta'$ depending only on $L,\vartheta$ and $\theta$.
\end{lemma}
\bpf We shall reduce this situation to the quasisymmetric case by using auxiliary inversions as in \cite{Vai-5}. By \cite[Theorem 3.19]{Vai-0}, $f$
  extends  continuously to a quasim\"{o}bius homeomorphism $f:\, \overline{D}\to \overline{D'}$. Without loss of generality with auxiliary translations if necessary, we may assume that $0\in \partial D$, and that $f(0)$ is either $0$ or $\infty$.

Let $u:\,\dot{E}\to \dot{E}$ be the inversion $u(x)=x/|x|^2$ for all $x\in\dot{E}$. Similarly,
 let $u':\,\dot{E'}\to \dot{E'}$ denote the inversion $u'(x')=x'/|x'|^2$ for all $x'\in \dot{E'}$. If $f(0)=0$, then we define $g:\,u(\overline{D})\to u'(\overline{D'})$ by $g=u'\circ f\circ u^{-1}$. If $f(0)=\infty$, then we define $g:\,u(\overline{D})\to \overline{D'}$ by $g=f\circ u^{-1}$. Thus, $g$ is $\theta_1$-quasim\"{o}bius with $\theta_1(t)=81\theta(81t)$ by \cite[Theorem 6.22]{Vai-5}. Since $g(\infty)=\infty$, $g$ is in fact $\theta_1$-quasisymmetric by means of \cite[Theorem 3.10]{Vai-0}.

On the other hand, by Lemma \ref{z-1.0} we note that the conclusion of the lemma is preserved under $L$-QH mappings. So, by \cite[5.11]{Vai-5}, no generality is lost in assuming that $f:\,\overline{D}\to \overline{D'}$ is $\theta$-quasisymmetric and $f:\,D\to D'$ is $(L,\vartheta)$-locally bilipschitz.

Let $\vartheta'\in(0,1)$ be a constant such that $2\eta(\vartheta')\leq \vartheta$.
Let $x',y'\in D'$ be two distinct points with $|x'-y'|< \vartheta' d_{D'}(x')$. Because $f$ is $(L,\vartheta)$-locally bilipschitz, it suffices to show
\be\label{z-1.8}
|x-y|<\vartheta d_D(x).
\ee
Choose $w\in\partial D$ with $|x-w|\leq 2d_D(x)$. Since $f^{-1}$ is $\eta$-quasisymmetric with $\eta(t)=\theta^{-1}(1/t)^{-1}$, we have
\[
\frac{|x-y|}{|x-w|}\leq \eta\Big(\frac{|x'-y'|}{|x'-w'|}\Big)\leq \eta(\vartheta'),
\]
and so $|x-y|\leq 2\eta(\vartheta')d_D(x)$. Furthermore, the inequality (\ref{z-1.8}) follows from the choice of $\vartheta'$ with $2\eta(\vartheta')\leq \vartheta$.
The proof of the lemma is complete.
\epf

\subsection{}
The final section is devoted to the proofs of Theorems \ref{main-thm-3}, \ref{main-thm-4} and \ref{main-thm-5}. To this end, we begin with an auxiliary lemma for our purpose. Throughout this subsection, we assume that $D\subsetneq E$ and $D'\subsetneq E'$ are proper domains, and
that $f:\,\overline{D}\to \overline{D'}$ is a homeomorphism. We also set $d(x)=d_D(x)$ and $d(x')=d_{D'}(x')$ for all $x\in D$ and $x'\in D'$.

\begin{lemma}\label{z-4}
Let $C,L\geq 1$ and $\theta:\,[0,\infty)\to [0,\infty)$ be a homeomorphism, and let $\partial D$ be $C$-uniformly perfect. If $f:\,\overline{D}\to \overline{D'}$ is $\theta$-quasim\"{o}bius relative to $\partial D$ and $f:\,\partial D\to \partial D'$ is $L$-bilipschitz, then there are positive constants $M_1$ and $M_2$ depending only on $\theta,C$ and $L$ such that
\begin{enumerate}
\item\label{z-5} $d(x)/M_1\leq d(x')\leq |x'-x_0'|\leq M_1d(x)$ for each $x\in D$ and $x_0\in\partial D\cap \overline{B}(x,2d(x))$,
\item\label{z-6} $|x-y|/M_2\leq|x'-y'|\leq M_2|x-y|$  for each $x\in D$ and $y\in\partial D$.
\end{enumerate}
\end{lemma}

\bpf We first prove a basic result by using the uniform perfectness of the boundary $\partial D$ as follows.
\bcl\label{zz-03} For each $x\in D$ and $x_0\in\partial D\cap \overline{B}(x,2d(x))$, there is an $x_1\in\partial D$ such that
\[ 
\frac{d(x)}{2C}\leq |x_0-x_1|\leq 6d(x).
\]
\ecl
If $\partial D\subset \overline{B}(x,3d(x))$, then we observe that $D$ is bounded. Choose a point $x_1\in\partial D$ with $|x_0-x_1|\geq \diam (\partial D)/2$. It follows from Lemma \ref{z-1} that
\[
|x_0-x_1|\geq \diam (\overline{D})/2\geq d(x)/2\;\;\mbox{and}\;\;|x_0-x_1|\leq |x_0-x|+|x-x_1|\leq 6d(x),
\]
as required.

Next we consider the remaining case that  $\partial D\setminus \overline{B}(x,3d(x))\not=\emptyset$. Since $x_0\in \overline{B}(x,2d(x))$, we have
\[
\overline{B}(x_0,d(x))\subset \overline{B}(x,3d(x))
\]
and so
\[
\partial D\setminus \overline{B}(x_0,d(x))\not=\emptyset.
\]
Moreover, since $\partial D$ is $C$-uniformly perfect, there is a point $x_1\in\partial D$ such that
\[
\frac{d(x)}{2C}\leq |x_0-x_1|\leq d(x)/2.
\]
Hence Claim \ref{zz-03} holds.

Next we verify (\ref{z-5}). We prove that there exists some $M_1$ such that $|x'-x_0'|\leq M_1d(x)$ for all $x\in D$, because the other direction follows similarly. By Corollary \ref{s-8}, we see that there is a homeomorphism $\theta_1:\,[0,\infty)\to [0,\infty)$ depending only on $\theta$ and $L$ such that $f:\,\overline{D}\to \overline{D'}$ is $\theta_1$-quasisymmetric relative to $\partial D$, because every $L$-bilipschitz map is $\eta_0$-quasisymmetric with $\eta_0(t)=L^2t$. Moreover, by Claim \ref{zz-03} we find that there is $x_1\in\partial D$ such that
\be\label{z-0.1}
\frac{d(x)}{2C}\leq |x_0-x_1|\leq 6d(x).
\ee
It follows that
\be\label{z-0.2}
\frac{|x'-x_0'|}{|x_1'-x_0'|}\leq \theta_1\Big(\frac{|x-x_0|}{|x_1-x_0|}\Big)\leq \theta_1(4C).
\ee
Now by (\ref{z-0.1}), (\ref{z-0.2}) and by the requirement that $f:\,\partial D\to \partial D'$ is $L$-bilipschitz, we obtain
\[
|x'-x_0'|\leq \theta_1(4C)|x_1'-x_0'|\leq L\theta_1(4C)|x_1-x_0|\leq 6L\theta_1(4C)d(x),
\]
as desired. This ensures (\ref{z-5}) with the choice of $M_1=6L\theta_1(4C)$.

Next we prove (\ref{z-6}). Fix $x\in D$ and $y\in \partial D$. We first show the right hand side inequality of (\ref{z-6}). If $y\in\partial D \cap \overline{B}(x,2d(x))$,
then by (\ref{z-5}) we see  that
\[
|x'-y'|\leq M_1d(x)\leq M_1|x-y|,
\]
as required.

So we may assume that $y\in\partial D$ with $|y-x|>2d(x)$. Let $0<\varepsilon<1/2$ and take a point $z\in\partial D$ with $|x-z|\leq (1+\varepsilon)d(x)$. A similar argument as above gives
\be\label{z-0.3}
|x'-z'|\leq M_1|x-z|.
\ee
Moreover, we have $|x-y|\geq|x-z|$ and so
\be\label{z-0.4}
|z-y|\leq |z-x|+|x-y|\leq 2|x-y|.
\ee
Since $f:\,\partial D\to \partial D'$ is $L$-bilipschitz, by (\ref{z-0.3}) and (\ref{z-0.4}), we get
\[
|x'-y'|\leq |x'-z'|+|z'-y'|\leq M_1|x-z|+L|z-y|\leq (M_1+2L)|x-y|,
\]
and we are done.

It remains to verify the left hand side inequality of (\ref{z-6}). On the one hand, one observes that the inverse of $f$ is $\theta_1'$-quasisymmtric relative to $\partial D'$ with $\theta_1'(t)=\theta_1^{-1}(1/t)^{-1}$. On the other hand, since $f:\,\partial D \to \partial D'$ is $L$-bilipschitz and  $\partial D$ is $C$-uniformly perfect, we see that $\partial D'$ is $(CL^2)$-uniformly perfect. A symmetric argument as above gives the left hand side inequality, from which (\ref{z-6}) follows by the choice of
\[
 M_2=\big(2L+6L\theta_1(4C)\big)\vee\big(2L+6L\theta_1'( 4CL^2)\big).
\]
The proof of the lemma is complete.
\epf

\subsection*{Proof of Theorem \ref{main-thm-4}:}
 Suppose that $\partial D$ is $C$-uniformly perfect, and that $f:\,\overline{D}\to \overline{D'}$ is a homeomorphism such that $f:\,D\to D'$ is $L$-QH, $f:\, \partial D\to \partial D'$ is $L$-bilipschitz, and $f:\, \overline{D}\to \overline{D'}$ is $\theta$-quasim\"obius relative to $\partial D$. We need to show that $f$ is bilipschitz on $\overline{D}$.

 First, by Lemma \ref{z-4}(\ref{z-6}) we only need to prove that $f:\, D\to D'$ is bilipschitz. It suffices to find a constant $L'>0$ such that
\[
|x'-y'|\leq L'|x-y| ~\mbox{ for all $x, y\in D$},
\]
because the other direction follows from a symmetric argument.

To this end, we divide the proof into two cases.

We first assume that $|x-y|< \mu d(x)$, where $\mu$ is the constant of Lemma \ref{z-3} depending only on $L$. Since $f:\,D\to D'$ is $L$-QH, it follows from Lemma \ref{z-3} that there is a constant $C>0$ such that
\[
\frac{|x'-y'|}{d(x')}\leq C\frac{|x-y|}{d(x)}.
\]
From this and Lemma \ref{z-4}(\ref{z-5}) it follows that there is a positive constant $M_1$ depending only on $\theta,C$ and $L$ such that
\[
|x'-y'|\leq CM_1|x-y|,
\]
as desired.

We are thus left to consider the case that $|x-y|\geq \mu d(x)$. Choose a point $x_0\in\partial D$ with $|x-x_0|\leq 2d(x)$. Then by Lemma \ref{z-4}(\ref{z-6}) we see that there is a positive constant $M_2$ depending only on $\theta,C$ and $L$ such that
\[
|x'-y'|\leq |x'-x_0'|+|x_0'-y'|\leq M_2(|x-x_0|+|x_0-y|)\leq M_2(4/\mu+1)|x-y|,
\]
where the last inequality follows from
\[
|x-x_0|+|x_0-y|\leq 2|x-x_0|+|x-y|\leq 4d(x)+|x-y|\leq (4/\mu+1)|x-y|.
\]
The proof of the theorem is  complete .
\qed

\vspace{6pt}

\subsection*{Proof of Theorem \ref{main-thm-5}.}
 Suppose that $\partial D$ is $C$-uniformly perfect, and   $f:\,\overline{D}\to \overline{D'}$ is a homeomorphism such that $f:\,D\to D'$ is $\theta$-quasim\"obius and  $(L, \vartheta)$-locally bilipschitz, and $f:\, \partial{D}\to \partial{D'}$ is $L$-bilipschitz.

First, by Lemma \ref{z-2.0} we deduce that $f^{-1}$ is $(L',\vartheta')$-locally bilipschitz for some constants $L',\vartheta'$ depending only on $\theta$, $L$ and $\vartheta$.
This, together with Lemma \ref{z-3}, guarantees that $f$ is $M$-QH, where $M$ depends only on $\theta$, $L$ and $\vartheta$.
Theorem \ref{main-thm-5} follows from this fact and Theorem \ref{main-thm-4}.
 \qed

\vspace{8pt}
Before the proof of Theorem \ref{main-thm-3}, we also need the following result.

\begin{lemma}\label{zz-04} Let $C,L\geq 1$, $\vartheta\in(0,1)$ and  $\varphi:\,[0,\infty)\to [0,\infty)$ be a homeomorphism. Let $D\subset E$ and $D'\subset E'$ be domains. Suppose that
\begin{enumerate}
\item $f:\,D\to D'$ is $\varphi$-FQC  and  $f:\,D\to D'$ is $(L,\vartheta)$-locally bilipschitz,
\item $f:\,\partial D\to \partial D'$ is $L$-bilipschitz,
\item $D$ is $c$-uniform and $\partial D$ is $C$-uniformly perfect.
\end{enumerate}
Then $f:\,\overline{D}\to \overline{D'}$ is $L'$-bilipschitz with $L'$ depending only on $c,C,L,\vartheta$ and $\varphi$.
\end{lemma}

\bpf First, by Theorem \ref{main-thm-1} we note that $D'$ is $c'$-uniform with $c'=c'(c,L,\varphi)$ because each $\varphi$-FQC mapping is $(M,C)$-CQH by \cite[Theorem 4.7]{Vai-5}. Then we deduce from Theorem \ref{Thm-A} that $f$ is $\theta$-quasim\"obius with $\theta=\theta(c,L,\varphi)$. Moreover, from this fact and Lemma \ref{z-2.0} we conclude that $f^{-1}$ is $(L',\vartheta')$-locally bilipschitz for some constants $L',\vartheta'$ depending only on $\theta$, $L$ and $\vartheta$. This, together with Lemma \ref{z-3}, shows that $f$ is $M$-QH, where $M$ depends only on $\theta$, $L$ and $\vartheta$. Furthermore, we use Theorem \ref{main-thm-4} to obtain that $f$ is $L'$-bilipschitz with $L'=L'(c,L,C,\vartheta,\varphi)$.
\epf

\subsection*{Proof of Theorem \ref{main-thm-3}.}

On the one hand, since $f$ is $L$-QH, by Lemma \ref{z-3} we know that $f$ is $(L_0,\vartheta)$-locally bilipschitz with $L_0$ and $\vartheta$ depending only on $L$. On the other hand, from \cite[Corollary 5.13]{Vai-5} it follows that $f$ is $\varphi$-FQC with $\varphi(t)=4L^2t$, because $f$ is $L$-QH.
Now, Theorem \ref{main-thm-3} follows from these two facts and Lemma \ref{zz-04}.
\qed

\vspace{6pt}

Finally, we consider Question \ref{q-3} for domains in the extended spaces and our result is as follows.

\begin{theorem}\label{zz-05} Let $C,L\geq 1$, and   $\theta:\,[0,\infty)\to [0,\infty)$ be a homeomorphism. Let $D\subset \dot{E}$ and $D'\subset \dot{E'}$ be domains. Suppose that
\begin{enumerate}
\item  $f:\,D\to D'$ is  $L$-QH and  $f:\,\overline{D}\to \overline{D'}$ is $\theta$-quasim\"{o}bius relative to $\partial D$,
\item $f:\,\partial D\to \partial D'$ is  $L$-bilipschitz and  $\partial D$ is $C$-uniformly perfect with $\diam (\partial D)> 0$,
\item\label{s-3} $\infty\in D$ if and only if $f(\infty)=\infty$ if and only if $\infty\in D'$.
\end{enumerate}
Then $f:\,\overline{D}\to \overline{D'}$ is $L'$-bilipschitz with $L'$ depending only on $C,L,\theta$.
\end{theorem}

\bpf We divide the proof into two cases.

 First, we let $\infty \not\in D$. Then $D$ is a domain of $E$. By (\ref{s-3}), we have $\infty \not\in D'$  which implies that  $D'$ is also a domain of $E'$. The conclusion follows from Theorem \ref{main-thm-4} in this case.

 Secondly, we let $\infty \in D$. Then, again by (\ref{s-3}), we get $f(\infty)=\infty \in D'$. By auxiliary translations, we may assume that $0\in\partial D$ and $f(0)=0\in\partial D'$. Let $u:\,\dot{E}\to \dot{E}$ be the inversion $u(x)=x/|x|^2$ for all $x\in \dot{E}$. Similarly, let $u':\,\dot{E'}\to \dot{E'}$ be the inversion $u(x')=x'/|x'|^2$ for all $x'\in \dot{E'}$. Then it follows from \cite[(5.9), 5.11, 5.12]{Vai-5} that
both $u_1=u|_D$ and $u_2=u'|_{D'}$ are $\theta_0$-quasim\"{o}bius and $36$-QH with $\theta_0(t)=81t$.  Next, we let $G=u_1(D)$, $G'=u_2(D')$ and
\[
g=u_2\circ f\circ u_1^{-1}:\,G\to G'.
\]

On the one hand, we see from the above facts that $g:\,G\to G'$ is $L_1$-QH with $L_1=1296L$ and $\theta_1$-quasim\"{o}bius with $\theta_1(t)=81\theta(81t)$, because the composition of QH (resp. quasim\"obius) maps  is QH (resp. quasim\"obius).

On the other hand, we observe that $g$ is $L_2$-bilipschitz on $\partial G$ with $L_2=9L^3$. Indeed,  using \cite[(5.9)]{Vai-5} we compute
\[
|g(p)-g(q)|\leq \frac{3|f\circ u_1^{-1}(p)-f\circ u_1^{-1}(q)|}{|f\circ u_1^{-1}(p)||f\circ u_1^{-1}(q)|}\leq \frac{3L^3|u_1^{-1}(p)-u_1^{-1}(q)|}{|u_1^{-1}(p)|| u_1^{-1}(q)|}\leq 9L^3|p-q|
\]
for all $p,q\in\partial G$ with $\infty\not\in\{p,q\}$,  because $0=f(0)\in f\circ u_1^{-1}(\partial G)=\partial D'.$

Moreover, we have $\infty\in\partial G$ and $g(\infty)=\infty \in \partial G'$. This implies that $G\subset E$ and $G'\subset E'$. Furthermore, we find that $\partial G=u_1(\partial D)$ is $C_1$-uniformly perfect with $C_1$ depending only on $C$, because quasim\"{o}bius maps preserve the uniform perfectness by \cite[Theorem C]{WZ17}. Hence we may apply Theorem \ref{main-thm-4} to obtain that $g:\,\overline{G}\to \overline{G'}$ is $L_3$-bilipschitz.

It remains to show that $f=u_2^{-1}\circ g\circ u_1$ is $L'$-bilipschitz with $L'=9L_3^3$. Again, by \cite[(5.9)]{Vai-5}, for all $x,y\in \overline{D}$ with $\infty\not \in \{x,y\}$, we have
\[
|f(x)-f(y)|\leq \frac{3|g\circ u_1(x)-g\circ u_1(y)|}{|g\circ u_1(x)||g\circ u_1(y)|}\leq \frac{3L_3^3|u_1(x)-u_1(y)|}{|u_1(x)|||u_1(y)}\leq 9L_3^3|x-y|,
\]
because $g(0)=u_2\circ f(\infty)=u_2(\infty)=0\in G'$.
Hence we are done and the proof of this theorem is complete.
\epf

We conclude this section by an example which infers that the condition (\ref{s-3}) of Theorem \ref{zz-05} is necessary.

\beg\label{exa-3}
Let $\mathbb{\dot{C}}=\mathbb{C}\cup \{\infty\}$ and $S_0=\{e^{it}\in \mathbb{\dot{C}}:\, -\pi/2\leq t\leq \pi/2\}$, and let $D=D'=\mathbb{\dot{C}}\setminus\{S_0\}$. Denote $u(z)=z^{-1}$ for $z\in \dot{\mathbb{C}}$.
Thus $u|_{S_0}=\id_{S_0}$, $u$ is a M\"{o}bius and $L$-QH mapping. However, it is clear that $u$ is not $L$-bilipschitz for any $L<\infty$.
\eeg

\bigskip

\subsection*{Acknowledgements} The authors are indebted to the referee for the valuable suggestions.
Qingshan Zhou was partly supported by NNSF of
China (Nos. 11901090 and 12071121), by Department of Education of Guangdong Province, China (No. 2021KTSCX116), by Guangdong Basic and Applied Basic Research Foundation (No. 2021A1515012289).


\normalsize


\begin{thebibliography}{[HD82]}




\normalsize
\baselineskip=17pt




\bibitem{BB03} {\rm Z. M. Balogh and S. M. Buckley}, \emph{Geometric characterizations of Gromov hyperbolicity,}\textrm{Invent. Math.}, {\rm 153} (2003), 261--301.

\bibitem{BHK}  {\rm M. Bonk, J. Heinonen, and P. Koskela},
\emph{Uniformizing Gromov hyperbolic spaces,} \textrm{Ast\'{e}risque,} {\rm 270} (2001), viii+99 pp.

\bibitem{BK1}  {\rm M. Bonk and B. Kleiner},
\emph{Quasisymmetric parametrizations of two-dimensional metric spheres,} 
\textrm{Invent. Math.,} {\rm 150} (2002),  127--183.

\bibitem{Bu}  {\rm C. Butler}, \emph{Uniformization, $\partial$-biLipschitz maps, sphericalization, and inversion,}
\textrm{arXiv:2008.06806.}


\bibitem{GO}  {\rm F. W. Gehring and B. G. Osgood}, 
\emph{Uniform domains and the quasi-hyperbolic metric,} 
 \textrm{J. Anal. Math.,} {\rm 36} (1979),
50--74.

\bibitem{GP}  {\rm F. W. Gehring and B. P. Palka}, 
\emph{Quasiconformally homogeneous domains,} 
\textrm{J. Anal. Math.,} {\rm 30} (1976), 172--199.

\bibitem{Gr87} {\rm M. Gromov}, Hyperbolic groups, \emph{Essays in group theory}, Math. Sci. Res. Inst. Publ., Springer, pp. 75--263, 1987.



\bibitem{Her06}  {\rm D. Herron},  \emph{Quasiconformal deformations and volume growth,} 
\textrm{Proc. London Math. Soc.}, {\rm 92} (2006), 161--199.

\bibitem{HSX}  {\rm D. Herron, N. Shanmugalingam, and X. Xie}, 
\emph{Uniformity from Gromov hyperbolicity,} \textrm{Illinois J. Math.}, {\rm 52} (2008), 1065--1109.


\bibitem{HL}  {\rm X. Huang and J. Liu}, 
\emph{Quasihyperbolic metric and quasisymmetric mappings in metric spaces,} 
{\it Trans. Amer. Math. Soc.,} {\rm 367} (2015), 6225--6246.

\bibitem{HLL}  {\rm X. Huang, H. Liu, and J. Liu},
\emph{Local properties of quasihyperbolic mappings in metric spaces,}
 \textrm{ Ann. Acad. Sci. Fenn. Ser. A I Math.,} {\rm 41} (2016), 23--40.

\bibitem{HLVW}  {\rm M. Huang, Y. Li, M. Vuorinen, and X. Wang}, 
\emph{On quasim\"obius maps in real Banach spaces,} \textrm{Israel J. Math.,} {\rm 198} (2013), 467--486.

\bibitem{HRWZ}  {\rm M. Huang, A. Rasila, X. Wang, and Q. Zhou}, 
\emph{Semisolidity and locally weak quasisymmetry of homeomorphisms in metric spaces,}
 \textrm{Studia Math.,} {\rm 242} (2018), 267--301.

\bibitem{K}  {\rm B. Kleiner}, 
\emph{The asymptotic geometry of negatively curved spaces: uniformization, geometrization and rigidity,}
\textrm{International Congress of Mathematicians. Vol. II, 743--768, Eur. Math. Soc., Z\"urich,} (2006).



\bibitem{LVX}{\rm Y. Li, M. Vuorinen, and X. Wang}, 
\emph{Quasiconformal maps with bilipschitz or identity boundary values in Banach spaces.}
\textrm{Ann. Acad. Sci. Fenn. Math.,} {\rm 39} (2014), 905--917.

\bibitem{rt} {\rm A. Rasila and J. Talponen,} 
\emph{Convexity properties of quasihyperbolic balls on Banach spaces,} 
Ann. Acad. Sci. Fenn.,  {\rm 37} (2012), 215--228.

\bibitem{rt2} {\rm A. Rasila and J. Talponen,}
\emph{On Quasihyperbolic Geodesics in Banach Spaces,}  \textrm{Ann. Acad. Sci. Fenn. Ser. A I Math.}, {\rm 39} (2014), 163--173.


\bibitem{TV}  {\rm P. Tukia and J. V\"{a}is\"{a}l\"{a}}, 
\emph{Quasisymmetric embeddings of metric spaces,}
\textrm{Ann. Acad. Sci. Fenn. Ser. A I Math.,} {\rm 5} (1980), 97--114.

\bibitem{Vai-0}  {\rm J. V\"{a}is\"{a}l\"{a}}, 
\emph{Quasim\"obius maps,} 
\textrm{J. Anal. Math.,} {\rm 44} (1985), 218--234.

\bibitem{Va90}  {\rm J. V\"{a}is\"{a}l\"{a}}, 
\emph{Quasisymmetry and unions,} 
\textrm{Manuscripta Math.,} {\rm 68} (1990), 101--111.

\bibitem{Vai-1}  {\rm J. V\"{a}is\"{a}l\"{a}}, 
\emph{Free quasiconformality in Banach spaces. I,} 
\textrm{Ann. Acad. Sci. Fenn. Ser. A I Math.,} {\rm 15} (1990), 355--379.

\bibitem{Vai-2}  {\rm J. V\"{a}is\"{a}l\"{a}}, 
\emph{ Free quasiconformality in Banach spaces. II,} 
\textrm{Ann. Acad. Sci. Fenn. Ser. A I Math.,} {\rm 16} (1991), 255--310.

\bibitem{Vai-3}  {\rm J. V\"{a}is\"{a}l\"{a}}, \emph{Free quasiconformality in Banach spaces. III,} 
\textrm{Ann. Acad. Sci. Fenn. Ser. A I Math.,}
{\rm 17} (1992), 393--408.

\bibitem{Vai-4}  {\rm J. V\"{a}is\"{a}l\"{a}}, 
\emph{Free quasiconformality in Banach spaces, IV,} 
\textrm{Analysis and Topology, 697--717, World Sci. Publ., River Edge, N. J.,} (1998).


\bibitem{Vai-5} {\rm J. V\"{a}is\"{a}l\"{a}}, 
\emph{The free quasiworld. Freely quasiconformal and related maps in Banach spaces,}
\textrm{Quasiconformal geometry and dynamics $($Lublin 1996$)$, Banach Center Publications, Vol. 48,} (1999), 55--118.


\bibitem{Vai}  {\rm J. V\"{a}is\"{a}l\"{a}}, 
\emph{Hyperbolic and uniform domains in Banach spaces,}
\textrm{Ann. Acad. Sci. Fenn. Ser. A I Math.,} {\rm 30} (2005), 261--302.

\bibitem{Vai-00}  {\rm J. V\"{a}is\"{a}l\"{a}}, \emph{Gromov hyperbolic spaces,}
\textrm{Expo. Math.}, {\rm 23} (2005), 187--231.

\bibitem{V16} {\rm V. Vellis}, 
\emph{Extension properties of planar uniform domains,}
\textrm{arXiv:1609.08763} (2016).




\bibitem{WZ17}  {\rm X. Wang and Q. Zhou}, 
\emph{Quasim\"{o}bius maps, weakly quasim\"{o}bius maps and uniform perfectness in quasi-metric spaces,} 
\textrm{Ann. Acad. Sci. Fenn. Math.}, {\rm 42} (2017), 257--284.

\bibitem{Z}  {\rm Q. Zhou}, 
\emph{Quasihyperbolic mappings in Banach spaces,} 
\textrm{Ann. Fenn. Math.}, {\rm 46} (2021), 335--344.

\bibitem{ZR}  {\rm Q. Zhou and A. Rasila}, 
\emph{Quasim\"obius invariance of uniform domains,} 
\textrm{Studia Math.}, {\rm 261} (2021), 1--24.




\end{thebibliography}
\end{document}